\documentclass[11pt]{amsart}
\usepackage{amsmath}
\usepackage{amssymb}
\usepackage{amscd}

\def\NZQ{\mathbb}               
\def\NN{{\NZQ N}}
\def\QQ{{\NZQ Q}}
\def\ZZ{{\NZQ Z}}
\def\RR{{\NZQ R}}
\def\CC{{\NZQ C}}
\def\AA{{\NZQ A}}

\newcommand{\bfind}[1]{\index{#1}{\bf #1}}

\newcommand{\sn}{\par\smallskip\noindent}

\newcommand{\pars}{\par\smallskip}
\newcommand{\parm}{\par\medskip}

\newcommand{\sep}{^{\it sep}}
\newcommand{\chara}{\mbox{\rm char}\,}

\newcommand{\cO}{\mathcal{O}}

\newcommand{\Z}{\mathbb Z}

\newtheorem{theorem}{Theorem}[section]

\theoremstyle{definition}

\newtheorem{Theorem}{Theorem}[section]
\newtheorem{Lemma}[Theorem]{Lemma}
\newtheorem{Corollary}[Theorem]{Corollary}
\newtheorem{Proposition}[Theorem]{Proposition}

\newtheorem{Example}[Theorem]{Example}

\newtheorem{Definition}[Theorem]{Definition}

\newtheorem{Question}[Theorem]{Question}
\let\epsilon\varepsilon
\let\phi=\varphi
\let\kappa=\varkappa
\begin{document}

\title{A general theorem on the analytic spread and associated primes of divisorial filtrations}
\author{Steven Dale Cutkosky}
\thanks{Partially supported by NSF grant DMS-2348849}

\address{Steven Dale Cutkosky, Department of Mathematics,
University of Missouri, Columbia, MO 65211, USA}
\email{cutkoskys@missouri.edu}


\begin{abstract} 
We prove that a classical theorem of McAdam about the  analytic spread of an ideal in a Noetherian local ring $R$ is true for rational divisorial filtrations on an
arbitrary Noetherian local ring. We deduce McAdam's original theorem from our result.
This generalizes earlier work by the author, showing that this theorem is true for rational divisorial filtrations, with the restriction that $R$
 is an equidimensional local ring which is essentially of finite type over an excellent local ring of dimension $\le 3$.
\end{abstract}

\maketitle

\section{Introduction}

Notation used in this paper is explained in Section \ref{SecNot}.

McAdam proved the following important theorem in the theory of ideals. 
\begin{Theorem}\label{Theorem2}(\cite{McA}, Theorem 5.4.6 \cite{HS})  Let $R$ be a formally equidimensional local ring and $I$ be an ideal in $R$. Then $m_R\in \mbox{Ass}(R/\overline{I^n})$ if and only if the analytic spread $\ell(I)=\dim R$.
\end{Theorem}

In this paper, we generalize this theorem to divisorial filtrations in arbitrary Noetherian local rings. These filtrations are generally not Noetherian. Theorem \ref{Theorem1} is proven in Section \ref{SecThm1}.

 \begin{Theorem}\label{Theorem1} Let $R$ be a Noetherian local ring of dimension $d$, and $\mathcal I=\{I_n\}$ be a rational divisorial filtration of $R$. Then the following are equivalent. 
 \begin{enumerate}
 \item[1)] The analytic spread $\ell(\mathcal I)=d$.
 \item[2)] $m_R\in \mbox{Ass}_R(R/I_n)$ for some $n>0$.
 \item[3)] $m_R\in \mbox{Ass}_R(R/I_n)$ for all $n\gg 0$.
 \end{enumerate}
 \end{Theorem}
 
 Theorem \ref{Theorem2} is deduced from Theorem \ref{Theorem1} in Section \ref{SecThm2}.
 
 We reproduce a simple example from \cite{C1} showing that 
 the conclusion 2) implies 1) of Theorem \ref{Theorem1} do not hold for more general filtrations.

\begin{Example}(Example 1.5 \cite{C1}) There exist $\RR$-divisorial filtrations $\mathcal I$ such that 
the conclusions of Theorem \ref{Theorem1} are false. 
\end{Example}

\begin{proof} Let $R=k[[x]]$ be a power series ring in one variable over a field $k$. Let $I_n=(x^{\lceil n\pi\rceil})$ for $n\in \NN$, and $\mathcal I=\{I_n\}$. For fixed $n$, $r\lceil n\pi\rceil>\lceil rn\pi\rceil +1$ for some  $r\in \ZZ_{>0}$. Thus $f\in I_n$ implies $f^r\in m_RI_{nr}$. Thus $\ell(\mathcal I)=\dim R[\mathcal I]/m_RR[\mathcal I]=0<1=\dim R$.
\end{proof}

 In \cite{C1}, Theorem \ref{Theorem1} is proven with the assumption that $R$ is an excellent local domain of equicharacteristic zero or of dimension $\le 3$ (or more generally when resolution of singularities holds). In \cite{C2}, the theorem is proven more generally with the assumption that $R$ is an equidimensional local ring which is essentially of finite type over an excellent local ring of dimension $\le 3$.

In the papers \cite{C1} and \cite{C2}, geometric methods are used to realize a quotient of 
\linebreak $R[\mathcal I]/m_RR[\mathcal I]$ as a subalgebra of a finitely generated $k$-algebra where  $k=R/m_R$ (a subfinite $k$-algebra), such that this quotient contains a finitely generated $k$-algebra of dimension $d$. By a theorem of Kemp and Trung in \cite{KT}, this implies that $\ell(\mathcal I)\ge d$. Since $\ell(\mathcal I)\le d$ for all graded families of ideals on $R$ by a lemma in \cite{CPS}, this proves 2) implies 1), with the assumptions of those papers. The proof of 2) implies 1) in \cite{C1} used a relative form of the theory of $\sigma$-decomposition of pseudo effective divisors on a projective nonsingular variety of Nakayama, developed in \cite{N}. This theory requires resolution of singularities, which is not known to be true in dimension larger than three. This problem is finessed in \cite{C2}, by using a theorem of Temkin on existence of regular alterations \cite{T}, which allowed the theorem to be extended to rings which are essentially of finite type over an excellent local ring of dimension three. 

The proof in this paper also proves 2) implies 1) by showing that a quotient of 
\linebreak $R[\mathcal I]/m_RR[\mathcal I]$ is 
a subfinite $k$-algebra  which contains a finitely generated $k$-algebra of dimension $d$. What is different in this paper is that we show that a divisorial valuation $v$ which is essential for the filtration $\mathcal I$ is centered at $m_R$ and is such that the associated graded ring ${\rm gr}_v(R)$ is subfinite.  We then show that a quotient $C$ of $R[\mathcal I]/m_RR[\mathcal I]$ is a subalgebra of ${\rm gr}_v(R)$, and then 
use the fact that $\mbox{trdeg}_{k}k(v)=d-1$ to prove that $C$ contains a polynomial ring over $k$ of dimension $d$. 

The proofs of 1) implies 3) in \cite{C1} and \cite{C2} follow from Theorem 4.6 and the following remark in \cite{CPS}, and the generalization in Theorem 1.4 \cite{C2}.

\begin{Theorem}(\cite[Theorem 4.6]{CPS} and Theorem 1.4 \cite{C2}))\label{TheoremN2} Suppose that $R$ is a Noetherian local ring and $\mathcal I=\{I_n\}$ is  an $\RR$-divisorial filtration on $R$ such that $\ell(\mathcal I)=\dim R$. Then there exists $n_0\in \NN$ such that $m_R\in \mbox{Ass}(R/\overline{I^n})$ for all $n\ge n_0$.
 \end{Theorem} 

1) implies 3) of Theorem \ref{Theorem1} also follows from Theorem \ref{TheoremN2}. We give a different proof of 1) implies 3) of Theorem \ref{Theorem1} in this paper, at the end of Section \ref{SecThm1}.

Theorem \ref{TheoremN2} above generalizes Burch's theorem for ideals (\cite{Bu} and Proposition 5.4.7 \cite{HS}) to real divisorial filtrations.

\section{Notation and preliminaries}\label{SecNot} Let $R$ be a commutative ring. If $P$ is a prime ideal of $R$ then we denote the residue field of $P$ by
$\kappa(P)=(R/P)_P= \mbox{QF}(R/P)$, where $\mbox{QF}(A)$ denotes the quotient field of a domain $A$. If $R$ is a local ring, we denote the maximal ideal of $R$ by $m_R$. If $v$ is a valuation of a field $K$ we denote the valuation ring of $v$ by $\mathcal O_v$, its maximal ideal by $m_v$ and its residue field by $\kappa(v)=\mathcal O_v/m_v$. We denote the value group of $v$ by $\Gamma_v$.

\begin{Definition}(Definition 9.3.1 \cite{HS}) Let $R$ be a Noetherian integral domain with field of fractions $K$ and let $v$ be a valuation of $K$ such that $v$ is nonnegative on $R$. Let $P=m_v\cap R$. If $\mbox{trdeg}_{\kappa(P)}\kappa(v)=\mbox{ht }P-1$, then $v$ is said to be a divisorial valuation with respect to $R$.
\end{Definition}
The valuation ring of a divisorial valuation is Noetherian and has value group $\Gamma_v\cong \ZZ$, as shown in \cite{HS}. 

\begin{Definition} Suppose that $R$ is a $d$-dimensional  local ring. A divisorial valuation of $R$ is a divisorial valuation $v$ of the quotient ring $R/P$ where $P$ is a minimal prime of $R$ such that $\dim R/P=d$. We extend $v$ to a pseudo valuation  $v:R\rightarrow \Gamma_v\cup\{\infty\}$ where $v(f)=\infty$ if $f\in P$. 
\end{Definition}

Suppose that  $R$ is a  local ring and $v$ is a divisorial valuation of $R$, which we may assume has value group $\Gamma_v=\ZZ$. We define valuation ideals in $R$ for $\lambda\in \RR_{\ge 0}$ by
$I(v)_\lambda=\{f\in R\mid v(f)\ge \lambda\}$. $I(v)_1$ is a prime ideal $P$, called the center of $v$ on $R$. It is the preimage of $m_v$ in $R$. We also will say that $v$ is centered at $P$ on $R$.
$I(v)_0=R$, $I(v)_{\lambda}=I(v)_1$ if $0<\lambda\le 1$ and $I(v)_{\lambda}$ is $I(v)_1$-primary if $\lambda>0$.
We define the associated graded ring of the valuation by 
$$
{\rm gr}_v(R)=\oplus_{n\ge 0}I(v)_n/I(v)_{n+1}.
$$
The ring ${\rm gr}_v(R)$ is an integral domain, since $v(ab)=v(a)+v(b)$ for all $a,b\in R$.
For $f\in R$, the initial form of $f$ with respect to $v$ is the class ${\rm in}_v(f)$ of $f$ in 
$$
I(v)_{v(f)}/I(v)_{v(f)+1}\subset {\rm gr}_v(R).
$$

A graded family of ideals $\mathcal I=\{I_n\}_{n\in \ZZ_{\ge 0}}$ in a ring $R$ is a family of ideals $I_n$ in $R$ such that $I_mI_n\subset I_{m+n}$ if $m,n\in \ZZ_{\ge 0}$. The Rees algebra of a graded family of ideals $\mathcal I$ is $R[\mathcal I]=\oplus_{n\ge 0}I_n$ and the analytic spread of $\mathcal I$ is 
$\ell(\mathcal I)=\dim R[\mathcal I]/m_RR[\mathcal I]$. In the case that $\mathcal I$ is an $I$-adic filtration $\mathcal I=\{I^n\}$ where $I$ is an ideal of $R$, these definitions give the Rees algebra $R[I]$ and the analytic spread $\ell(I)$ of the ideal $I$. In a $d$-dimensional Noetherian local ring $R$, and graded family of ideals $\mathcal I$, we have that 
$\ell(\mathcal I)\le d$ and all values from 0 to $d$ can be realized, as shown in \cite{CPS}.

 Suppose that $R$ is a Noetherian local ring. A graded family of ideals $\mathcal I$ in $R$ is called a divisorial filtration if there exist $r\in\ZZ_{>0}$, $a_1,\ldots, a_r$ in $\RR_{>0}$ and divisorial valuations $v_1\ldots,v_r$ of $R$ such that
\begin{equation}\label{eq1}
I_n=\{f\in R\mid v_i(f)\ge n a_i\mbox{ for }1\le i\le r\}.
\end{equation}
The expression (\ref{eq1}) is called a presentation of the divisorial filtration $\mathcal I$. In general, there are infinity many different representations of $\mathcal I$. If a representation has the property that for each $j$ with $1\le j\le r$ the ideal $\{f\in R\mid v_i(f)\ge n a_i\mbox{ for }i\ne j\}$ properly contains $I_n$, then the representation is said to be irredundant. 

$\mathcal I$ is said to be a rational divisorial filtration if there exists a   presentation (\ref{eq1})  such that $a_i\in \QQ_{>0}$ for all $i$. $\mathcal I$ is an integral divisorial filtration if all $a_i\in \ZZ$.

\section{Proof of Theorem \ref{Theorem1}}\label{SecThm1} 

\begin{Lemma}\label{Lemma1} Let $S$ be a Noetherian local domain of dimension $e$ and let $v$ be a divisorial valuation of $S$ whose center is $m_S$.
Identify the value group of $v$ with $\ZZ$. Then there exists a complete local domain $T$ which is a quotient of the $m_S$-adic completion $\hat S$ and a divisorial valuation $w$ of $T$ such that $w| S=v$ and ${\rm gr}_v(S)\cong {\rm gr}_w(T)$ as graded $k=\kappa(m_R)$-algebras.
\end{Lemma}

\begin{proof} Let $V=\mathcal O_v$, which is a rank 1 discrete valuation ring. Let $\hat V$ be the completion of $V$ at its maximal ideal, with natural valuation $\hat v$ extending $v$. Since the center of $v$ on $S$ is $m_S$, we have that $v(m_S)\ge 1$, so that the natural homomorphism $S\rightarrow V$ is continuous for the $m_S$-adic and $m_V$-adic topologies, and thus extends to a homomorphism $\hat\Phi:\hat S\rightarrow \hat V$. Let $H$ be the kernel of $\Phi$ and $T=\hat S/H$. Then $H$ is a prime ideal, $H\cap S=0$ and $T$ is a complete Noetherian local domain embedded in $\hat V$. Let $w$ be the restriction of $\hat v$ to $\mbox{QF}(T)$. $w$ is nonnegative on $T$ and is centered at the maximal ideal of $T$. We have that 
$\kappa(m_T)=\kappa(m_R)$, $\Gamma_w=\Gamma_v=\ZZ$ and $\kappa(w)=\kappa(v)$. By Abhyankar's inequality, \cite{Ab} or Theorem 6.6.7 \cite{HS},
$$
e=1+{\rm trdeg}_k\kappa(v)\le \dim T\le \dim \hat S=e.
$$
Thus $\dim T=e$ and so $w$ is divisorial on $T$. The inclusions $S\rightarrow T$ give injective maps
\begin{equation}\label{eq2}
I(v)_j/I(v)_{j+1}\rightarrow I(w)_j/I(w)_{j+1}
\end{equation}
For $y\in I(w)_j$, there exists $x\in S$ such that $x=y+h$ with $h\in m_T^n$, with $n$ sufficiently large that $nw(m_T)>j$. Then $v(x)=j$ and $y-x\in I(w)_{j+1}$, so that (\ref{eq2}) is an isomorphism, giving the isomorphism ${\rm gr}_v(S)\cong {\rm gr}_w(T)$.
\end{proof}

Let $k$ be a field. A $k$-algebra is subfinite if it is a subalgebra of a finitely generated $k$-algebra.

\begin{Lemma}\label{Lemma2}
Let $R$ be an excellent local domain of dimension $d$ and $v$ be a divisorial valuation of the quotient field $K={\rm QF}(R)$ of $R$ whose center on $R$ is $m_R$. Let $k=\kappa(m_R)$. Then ${\rm gr}_v(R)$ is a subfinite $k$-algebra.
\end{Lemma}

\begin{proof} Let $\pi:X\rightarrow \mbox{Spec}(R)$ be such that $X$ is normal and there is a prime divisor $E\subset X$ such that $\mathcal O_{X,E}=\mathcal O_v$. Such an $X$ exists by Lemma 3.3 \cite{C2}. Since the center of $v$ on $R$ is $m_R$, $\pi(E)=m_R$ so the closed integral subscheme 
$E$ of $X$ is projective over $k$. 
Choose $g\in K^{\times}$ with $v(g)=1$ and set $B={\rm div}_X(g)-E$, where ${\rm div}_X(g)$ is the Weil divisor on $X$ associated to $g$, so the coefficient of $E$ in  $B$ is zero. Write $B=B_+-B_-$ where $B_+$ is the divisor of zeros of $B$ and $B_-$ is the divisor of poles. Since $X$ is projective over $\mbox{Spec}(R)$, there exists an effective ample Cartier divisor $H$ on $X$ such that $H\ge B_+$ and $E$ is not contained in the support of $H$. 

Suppose that  $f\in I(v)_j$. Then  $f$ is regular everywhere on $X$ and the order of vanishing of $f$ along $E$ is at least $j$. Thus ${\rm div}_X(f)\ge jE$, and so
$$
{\rm div}_X(\frac{f}{g^j})+jH={\rm div}_X(f)-j{\rm div}_X(g)+jH\ge j(H-B)\ge 0.
$$
Since $X$ is normal, this shows that
$$
\frac{f}{g^j}\in H^0(X,\mathcal O_X(jH)).
$$
Mapping $f$ to the image of $\frac{f}{g^j}$ after restriction to $E$, defines a map 
\begin{equation}\label{eq3}
I(v)_j\rightarrow H^0(E,\mathcal O_X(jH)\otimes_{\mathcal O_X}\mathcal O_E).
\end{equation}
$f\in I(v)_j$ is in the kernel of the map if and only if $v(\frac{f}{g^j})>0$ which holds if and only if $f\in I(v)_{j+1}$. Thus the kernel of (\ref{eq3}) is $I(v)_{j+1}$. These maps preserve multiplication. Thus they give an injective $k$-algebra homomorphism
$$
{\rm gr}_v(R)\rightarrow S=\oplus_{j\ge 0}H^0(E,\mathcal O_X(jH)\otimes_{\mathcal O_X}\mathcal O_E).
$$
The invertible sheaf $\mathcal O_X(jH)\otimes_{\mathcal O_X}\mathcal O_E$ is ample on the projective $k$-variety $E$. Hence $S$ is a finitely generated $k$-algebra proving Lemma \ref{Lemma2}.
\end{proof}

We now prove Theorem \ref{Theorem1}. We first prove that 2) implies 1). We can choose a presentation (\ref{eq1}) of $\mathcal I$ which is irredundant. Let $\mathfrak p_i$ be the center of $v_i$ on $R$, for $1\le i\le r$. Each ideal $I(v_i)_n$ is $\mathfrak p_i$-primary. The natural injection $R/I_n\rightarrow \oplus_iR/I(v_i)_n$ shows that $\mbox{Ass}_R(R/I_n)\subset \{\mathfrak p_1,\ldots,\mathfrak p_r\}$. By our assumption that 2) holds, we conclude that at least one valuation in the irredundant presentation is centered at $m_R$. 
Reindex the $v_i$ so that this valuation is $v_1$ and let $v=v_1$. Let 
$P=\{x\in R\mid v(x)=\infty\}$ and let $S=R/P$, so that $S$ is a $d$-dimensional local domain and $v$ is a divisorial valuation of $S$ centered at $m_S$. Let $a$ be the coefficient of $v=v_1$ in the defining inequalities of (\ref{eq1}). Let
$$
C=\sum_{n\ge 0}I_nt^n/(\sum_{n\ge 0}\{h\in I_n\mid v(h)>na\}t^n)=\sum_{n\ge 0}C_nt^n
$$
where $t$ is a variable of degree 1. We have a natural surjection of graded $k$-algebras
\begin{equation}\label{eq4}
R[\mathcal I]/m_RR[\mathcal I]\rightarrow C.
\end{equation}

If $na\in \ZZ$ then $I_n\cap I(v)_{na+1}=\{x\in I_n\mid v(x)>na\}$.
If $na\not\in \ZZ$ then every $h\in I_n$ has $v(h)\ge \lceil na \rceil>na$ so $C_n=0$. Thus we have a natural inclusion of $k$-algebras $C\rightarrow {\rm gr}_{v}S$. Thus $C$ is subfinite by Lemmas \ref{Lemma1} and \ref{Lemma2}. 

Since $v$ is divisorial and centered at $m_R$, we have that $\mbox{trdeg}_k\kappa(v)=d-1$. Choose algebraically independent elements 
$\xi_1,\ldots,\xi_{d-1}\in \kappa(v)$ over $k$. Let $\xi_0=1$. Lift $\xi_1,\ldots,\xi_d$ to elements of  $\mathcal O_v$, which necessarily have value zero, and lift $\xi_0$ to the 1 of $R$. After multiplying them by a suitable common element $h$ of $R$, we obtain $u_0,u_1,\ldots, u_{d-1}\in R$ which have the value
\begin{equation}\label{eq10}
v(u_0)=v(u_1)=\cdots=v(u_{d-1})=c>0,
\end{equation}
and the class of $\frac{u_i}{u_0}$ in $\kappa(v)$ is $\xi_i$ for $0\le i\le d-1$. Irredundancy of the $v_i$ implies that there exists an integer $N>0$ and $f\in R$ such that 
\begin{equation}\label{eq6}
v_i(f)\ge na_i\mbox{ for $i\ne 1$ and $v(f)<Na$}.
\end{equation} 
Set
\begin{equation}\label{eq7}
\delta=Na-v(f)\in \QQ_{>0}.
\end{equation}
Since $\delta$ is positive and rational, there exist positive integers $p$ and $q$ such that
\begin{equation}\label{eq12}
p\delta=qc.
\end{equation}
Let $h_i=f^pu_i^q$ for $0\le i\le d-1$. For every $v_i$ other than $v=v_1$, (\ref{eq6}) and nonnegativity on $R$ give
$$
v_j(h_i)=pv_j(f)+qv_j(u_i)\ge pNa_j.
$$
For $v$, equations  (\ref{eq10}), (\ref{eq7}) and (\ref{eq12}) give
$$
v(h_i)=pv(f)+qc=p(Na-\delta)+p\delta=pNa.
$$
Thus $h_i\in I_{pN}$ and each $h_it^{pN}$ has nonzero image $z_i$ in $C$. All the $z_i$ have the same positive degree $pN$.

We will now show that $z_0,\ldots,z_{d-1}$ are algebraically independent over $k$. Suppose otherwise, so that there exists a nonzero polynomial
$F(T_0,T_1,\ldots,T_{d-1})$ in the polynomial ring $k[T_0,\ldots,T_d]$ such that $F(z_0,\ldots,z_{d-1})=0$. Since $z_0,\ldots,z_{d-1}$ have the same degree, the homogeneous components of $F$ give homogeneous relations, so we may assume that $F$ is homogeneous, say of degree $b$. Now
$G={\rm gr}_v(S)$ is a domain. Since the natural map $C\rightarrow G$ is injective, $C$ is also a domain. Thus we have a natural inclusion
${\rm QF}(C)\rightarrow {\rm QF}(G)$ which takes $\frac{z_i}{z_0}$ to $\frac{{\rm in}_v(h_i)}{{\rm in}_v(h_0)}$.
Let 
$$
K_0(G)=\{\frac{a}{b}\in {\rm QF}(G)\mid a,b\mbox{ are nonzero homogeneous elements of $G$ of equal degree}\}\cup\{0\}.
$$
There is a canonical $k$-algebra isomorphism $K_0(G)\cong k(v)$ defined by 
$$
\frac{{\rm in}_v(f)}{{\rm in}_v(g)}\mapsto \overline{\left(\frac{f}{g}\right)},\mbox{ the class of 
$\frac{f}{g}$ in $\kappa(v)$},
$$
so we can identify $\frac{z_i}{z_0}$ with $\xi_i^q\in \kappa(v)$.
Since $F$ is homogeneous, $F(1,\frac{z_1}{z_0},\ldots,\frac{z_{d-1}}{z_0})=0$, so that $F(1,\xi_1^q,\ldots,\xi_{d-1}^q)=0$ in $\kappa(v)$. 
The algebraic independence of  $\xi_1^q,\ldots,\xi_{d-1}^q$ implies that $F(1,T_1,\ldots,T_{d-1})=0$. Since $F$ is homogeneous, we have that
$F(T_0,T_1,\ldots,T_{d-1})=0$ and so $F(T_0,T_1,\ldots,T_{d-1})=0$.
This contradiction shows that $z_0,z_1,\ldots,z_{d-1}$ are algebraically independent over $k$. Thus 
$$
k[T_0,\ldots,T_{d-1}]\cong k[z_0,\ldots,z_{d-1}]\subset C.
$$
Since $k[z_0,\ldots,z_{d-1}]\subset C$ are subfinite $k$-algebras, we have by Corollary 4.7 of \cite{KT} that 
$$
\dim C\ge \dim k[z_0,\ldots,z_{d-1}]=d.
$$ 
By (\ref{eq4}), we have that
$$
\ell_R(\mathcal I)=\dim R[\mathcal I]/m_RR[\mathcal I]\ge \dim C\ge d.
$$
Finally, $\ell_R(\mathcal I)\le \dim R=d$ by Lemma 3.6 \cite{CPS}, and so we have established that 2) implies 1).

We now prove that 1) implies 3). For $1\le i\le r$ and $n>0$, let
$$
\mathfrak p_i=\{h\in R\mid v_i(h)>0\}\mbox{ and }Q_{i,n}=\{h\in R\mid v_i(h)\ge na_i\}. 
$$
$\mathfrak p_i$ is the center of $v_i$ on $R$ and $Q_{i,n}$ is $\mathfrak p_i$-primary for $n>0$. 
Suppose that none of the $v_i$ have center $m_R$. By prime avoidance, there exists $x\in m_R\setminus\left((\cup_i\mathfrak p_i\right)\cup
\left(\cup_{P\in \mbox{Min}(R)}P\right)$, where the second union is over all minimal primes of $R$. Now $r\in (I_n:x)$ if and only if $v_i(xr)=v_i(r)\ge na_i$ for all $i$, so that 
$$(I_n:x)=I_n\mbox{ for all $n\ge 0$}.
$$
Thus 
\begin{equation}\label{eqL3}
I_n\cap xR=x(I_n:x)=xI_n.
\end{equation}
Let $R^*=R/xR$, $I_n^*=I_n(R/xR)$ and $\mathcal I^*=\{I_n^*\}$. By (\ref{eqL3}), there is a natural isomorphism
$$
R[\mathcal I]/xR[\mathcal I]\cong R^*[\mathcal I^*].
$$
Thus 
$$
R^*[\mathcal I^*]/m_RR[\mathcal I^*]\cong R[\mathcal I]/m_RR[\mathcal I]
$$
since $x\in m_R$. We have that $\dim R/xR\le d-1$ since $x$ avoids every minimal prime of $R$. Thus by Lemma 3.6 \cite{CPS}, we have that
$\ell(\mathcal I)=\ell(\mathcal I^*)\le d-1$, giving a contradiction. Thus some valuation, which we can assume is $v_1$, has center $\mathfrak p_1=m_R$.

Since our presentation of $\mathcal I$ is irredundant, there exists $N>0$ and $f\in R$ such that 
$$
v_i(f)\ge Na_i\mbox{ for $i>1$ and $c:=v_1(f)<Na_1$.}
$$
For each $n>0$, let $q_n=\lceil \frac{n}{N}\rceil$ and $h_n=f^{q_n}$. For every $i>1$,
$$
v_i(h_n)=q_nv_i(f)\ge q_nNa_i\ge na_i.
$$
Set $\epsilon=a_1-\frac{c}{N}>0$.  For $n\in\ZZ_{>0}$, 
$$
\left(\frac{n}{N}+1\right)c=n\frac{c}{N}+c=na_1-n\epsilon+c<na_1
$$
whenever $n>\frac{c}{\epsilon}$. Thus
$$
v_1(h_n)=q_nc< \left(\frac{n}{N}+1\right)c<na_1
$$
whenever $n>\frac{c}{\epsilon}$. Thus for all sufficiently large $n$,
$h_n\in\cap_{i\ne j}Q_{i,n}$, with $h_n\not\in Q_{1,n}$. Thus $h_n\not\in I_n$ and since $Q_{1,n}$ is $m_R$-primary, there exists $b_n>0$ such that $m_R^{b_n}\subset Q_{1,n}$ and so $m_R^{b_n}h_n\subset I_n$. Let $L=h_n(R/I_n)$. Let $s$ be the smallest power of $m_R$ such that $m_R^sL=0$. Let $0\ne z\in m_R^{s-1}L$. Then $m_Rz=0$ so $\mbox{Ann}(z)=m_R$. Thus $m_R\in \mbox{Ass}(R/I_n)$ for all $n$ sufficiently large, completing the proof of 1) implies 3).

The final implication 3) implies 2) is immediate, completing the proof of the theorem.

\section{Derivation of McAdam's Theorem \ref{Theorem2} from Theorem \ref{Theorem1}}\label{SecThm2}

Let $R$ be a formally equidimensional Noetherian local ring of dimension $d$ and $I$ be an ideal in $R$. Let $\mathcal I$ be the  filtration  $\mathcal I=\{\overline{I^n}\}$ of integral closures of powers of $I$.
Let 
$$
A=R[I]=\oplus_{n\ge 0}I^nt^n\mbox{ and }B=R[\mathcal I]=\oplus_{n\ge 0}\overline{I^n}t^n
$$
where $A$ and $B$ are graded by $\deg t=1$, be the Rees algebras of $I$ and $\mathcal I$ respectively. The natural inclusion $A\subset B$ is an integral extension. Thus the injection
$$
A/m_RB\cap A\rightarrow B/m_RB
$$
is also an integral extension. Let $\mathfrak p$ be a prime ideal of $A$ containing $m_RA$. Since $B$ is integral over $A$, there exists a prime ideal $\mathfrak q$ of $B$ such that $\mathfrak q\cap A=\mathfrak p$ by Theorem A.6 \cite{BH}. Thus $m_RB\cap A\subset \mathfrak q\cap A=\mathfrak p$, and so
$$
m_RB\cap A\subset \cap_{\mathfrak p \supset m_RA}\mathfrak p=\sqrt{m_RA}.
$$
Thus $\sqrt{m_RB\cap A}\subset \sqrt{m_RA}$. Certainly $\sqrt{m_RA}\subset \sqrt{m_RB\cap A}$, so that 
$$
\sqrt{m_RB\cap A}=\sqrt{m_RA}.
$$
Thus 
\begin{equation}\label{eqT1}
\ell(\mathcal I)=\dim B/m_RB=\dim A/m_RB\cap A=\dim A/m_RA=\ell(I),
\end{equation}
where the second equality is by Corollary A.8 \cite{BH}.

We first prove the theorem with the assumption that $I\not\subset P$ for every minimal prime ideal $P$ of $R$. 
The theory of Rees valuations of an ideal is developed in \cite{R} and Section 10 of \cite{HS}.
We have that $R/P$ is formally equidimensional of dimension $d$ (by Corollary B.4.3 \cite{HS}) and every localization of a formally equidimensional local ring ring is formally equidimensional (by Theorem B.5.2 \cite{HS}).  By Proposition 10.4.3 \cite{HS}, there exist $r>0$ and divisorial valuations $v_i$ on $R/P_i$ where $P_i$ is a minimal prime of $R$ and $a_i\in \ZZ_{>0}$ such that 
$$
\overline{I^n}=\{f\in R\mid v_i(f)\ge na_i\mbox{ for }1\le i\le r\}.
$$
It is implicitely assumed that $IV\ne 0$ for  $V=\mathcal O_v$ in this Proposition. 
This condition is satisfied since $I$ is not contained in a minimal prime ideal of $R$. The valuations $v_1,\ldots,v_r$ are called the Rees valuations of $I$.
The filtration $\mathcal I=\{\overline{I^n}\}$ is thus an integral divisorial filtration on $R$, so by Theorem \ref{Theorem1}, $m_R\in \mbox{Ass}(R/\overline{I^n})$ for some $n$ if and only if $\ell(\mathcal I)=d$. Since $\ell(\mathcal I)=\ell(I)$ by (\ref{eqT1}), Theorem \ref{Theorem2} is proven when $I$ is not contained in a minimal prime ideal of $R$.

Now assume that $I$ is any proper ideal of $R$ (possibly contained in a minimal prime of $R$). Let $R[X]$ be a polynomial ring over $R$ and let $S=R[X]_{\mathfrak n}$ where $\mathfrak n=(m_R,X)R[X]$. Then $S$ is a local ring with maximal ideal $m_S=\mathfrak n S$. Let $J=IS+XS$. The $m_S$-adic completion of $S$ is $\hat S\cong \hat R[[X]]$. Since $\hat R$ is equidimensional of dimension $d$, $\hat S$ is equidimensional of dimension $d+1$. $J$ contains the nonzero divisor $X$, so $J$ is not contained in any minimal prime of $S$. Thus 
\begin{equation}\label{eqT5}
m_S\in \mbox{Ass}_S(S/\overline{J^n})\mbox{ for some $n$ if and only if }\ell(J)=d+1
\end{equation}
by the first part of this proof. 

Let $\overline X$ be the class of $X$ in $F_J(S)=S[J]/m_SS[J]$.
$F_J(S)=\oplus_{n\ge 0}J^n/m_SJ^n$ where
$$
J^n/m_SJ^n\cong \oplus_{j=0}^n\left(I^{n-j}/m_RI^{n-j}\right)\overline X^j.
$$
Thus there is a natural ring isomorphism of the polynomial ring $F_I(R)[Z]$ over $F_I(R)=R[I]/m_RR[I]$ with $F_J(S)$, obtained by mapping $Z$ to $\overline X$. Hence
\begin{equation}\label{eqT6}
\ell(J)=\ell(I)+1.
\end{equation}
$$
\overline{J^n}=\overline{I^n}S+\overline{I^{n-1}}XS+\cdots+\overline IX^{n-1}S+X^nS.
$$
Thus as an $R$-module,
\begin{equation}\label{eqT8}
S/\overline{J^n}\cong \oplus_{j=0}^{n-1}(R/\overline{I^{n-j}})X^j.
\end{equation}
Let $M=S/\overline{J^n}$. Equation (\ref{eqT8}) implies that
\begin{equation}\label{eqT10}
L:=0:_M m_R\ne 0 \mbox{ if and only if }m_R\in \mbox{Ass}_R(R/\overline{I^s})\mbox{ for some }1\le s\le n.
\end{equation}
$X^jL\subset L$ for all $j\ge 0$ and $X^n\in J^n\subset \overline{J^n}$ since $X\in J$. Thus $X^nM=0$, and so $X^nL=0$.

Suppose that $L=0:_Mm_R\ne 0$ and let $z$ be a nonzero element of $L$.  Let $s\ge 1$ be the smallest positive integer such that $X^sz=0$. Then $w=X^{s-1}z\ne 0$. Thus $w\in 0:_Mm_R$ and $Xw=0$ so that $0\ne w\in 0:_Mm_S$ and thus $m_S\in\mbox{Ass}_S(M)=\mbox{Ass}_S(S/\overline{J^n})$.

Certainly $m_S\in \mbox{Ass}_S(M)$ implies $0:_Mm_R\ne 0$, so that 
$$
m_S\in \mbox{Ass}_S(M)\mbox{ if and only if }0:_Mm_R\ne 0.
$$
The conclusions of Theorem \ref{Theorem2} for arbitrary $I$ now follow from (\ref{eqT10}), (\ref{eqT5}) and (\ref{eqT6}).

\end{document}